\documentclass[]{article}
\usepackage{amsfonts, amsmath, amssymb}
\usepackage{graphicx,color}
\usepackage{graphics}
\usepackage{psfig}

\begin{document}

\title{A linear dispersive mechanism for numerical error growth: spurious caustics}

\author{Claire David \footnotemark[1] , Pierre Sagaut \footnotemark[1], Tapan Sengupta \footnotemark[2] \\
 \\
\footnotemark[1] Universit\'e Pierre et Marie Curie-Paris 6  \\
Laboratoire de Mod\'elisation en M\'ecanique, UMR CNRS 7607 \\
Bo\^ite courrier $n^0162$, 4 place Jussieu, 75252 Paris
cedex 05, France\\
 \\
\footnotemark[2] Indian Institute of Technology Kanpur  \\
Department of Aerospace Engineering \\
U.P. 208016, India }

\maketitle

\begin{abstract}

A linear dispersive mechanism for error focusing in polychromatic solutions is identified.
This local error pile-up corresponds to the existence of spurious caustics, which are allowed by the dispersive nature of the numerical error.
>From the mathematical point of view, spurious caustics are related to extrema of the  numerical group velocity.
Several popular schemes are analyzed and are shown to admit spurious caustics. It isalso observed that caustic-free schemes can be defined, like the Crank-Nicolson scheme.

\end{abstract}

\section{Introduction}
\label{sec:intro}

The analysis and control of numerical error in discretized propagation-type equations is of major
importance for both theoretical analysis and practical applications.
A huge amount of works has been devoted to the analysis of the numerical errors, its dynamics and 
its influence on the computed 
solution (the reader is referred to classical books, among which \cite{Hirsch,vich,lomax,tapan-book}).
It appears that existing works are mostly devoted to linear, one-dimensional numerical models, such as the linear advection equation

\begin{equation}
\label{transp}
 \frac{\partial u}{\partial t}+c \, \frac{\partial u}{\partial x}=0 
\end{equation}

where $c$ is a constant uniform advection velocity.
A striking observation is that, despite the tremendous efforts devoted to the analysis of 
numerical schemes in this simple case, the full exact non-homogeneous error equation has been
derived only very recently \cite{dipankar}.

The two sources of numerical error are the dispersive and dissipative properties of the numerical 
scheme, which are very often investigated thanks to the spectral
von Neumann analysis. Following this approach, a monochromatic wave is used to measure the accuracy of the scheme.
Such a tool is very powerful and provides the user with a deep insight into the discretization errors. 
But some results coming from practical numerical experiments still remain
unexplained, despite the linear character of the discrete numerical model. 
As an example, let us note the sudden growth of the numerical error for long range 
propagation reported by Zingg \cite{zingg} for a large set of numerical schemes, 
including optimized numerical schemes.

A limiting feature of the usual modal analysis is that it is applied to monochromatic reference 
solutions.
Therefore, dispersive phenomena associated to polychromatic solutions are usually not taken into 
account.

The present paper deals with the analysis of a linear dispersive mechanism which results in local  error focusing, i.e. to a sudden local error growth for polychromatic solutions.
This phenomena is reminiscent of the physical one referred to as the caustic phenomenon 
in linear dispersive physical models \cite{witham}, and will be referred to as the spurious caustic phenomenon hereafter. It will be shown that, for some specific values of the Courant-Friedrichs-Lewy (number (CFL), spurious caustics can exist for some popular finite-difference schemes.

The paper is organized as follows. 
The numerical schemes retained for the present analysis are briefly recalled in section 
\ref{sec:scheme}.
Main elements of cautic theory of interest for the present analysis are briefly recalled in 
section \ref{sec:caust}.
Test schemes are analyzed in section \ref{sec:analyze-test}.

\section{Test numerical schemes}
\label{sec:scheme}

For the sake of simplicity, the analysis will be restricted to three-level three-point numerical methods.
The extension of the present analysis to other schemes is straightforward. 
A general finite-differenced version of the linear advection equation (\ref{transp})  is

\begin{equation}
\label{eq:discrete-scheme}
\alpha  u_j^{n+1}+ \beta u_j^n + \gamma u_j ^{n-1}
+\delta u_{j+1}^n + \upsilon u_{j+1} ^{n-1}+ \varepsilon  u _{j-1} ^n + \zeta  u_{j+1}^{n+1}
+\eta u_{j-1}^{n-1} + \theta u_{j-1}^{n+1}  =0
\end{equation}
              
 with
    
\begin{equation}
{u_l}^m=u\,(l\,h, m\,\tau)
\end{equation}

where $h$ and  $\tau$ are the mesh size and time step respectively. For the sake of simplicity, these two quantities are assumed to be uniform.
The CFL number is defined as $\sigma = c \tau / h$, while the non-dimensional wave number is defined as $\varphi = k h$ where $k$ is the wave number of the signal under consideration.

A numerical scheme is  specified by selecting appropriate values of the coefficients  $\alpha$, $\beta$, $\gamma$, $\delta$,
$\varepsilon$, $\zeta$, $\eta$ and  $\theta$ in Eq. (\ref{eq:discrete-scheme}). Values corresponding to numerical schemes retained for the present works are given in Table 
\ref{tab:scheme}.

\begin{table}[htdp]
\caption{Numerical scheme coefficient.}
\begin{center}
\begin{tabular}{cccccccccc}
\hline
Name & $  \alpha $ & $ \beta$ & $\gamma$ & $\delta $  & $\epsilon$ & $\zeta$ & $\eta$ & $\theta$ & $\upsilon$ \\
\hline
Leapfrog & $\frac{1}{2 \tau} $ & 0 &  $\frac{-1}{2 \tau} $ & $  \frac{c}{ 2 h}  $ & $   \frac{-c}{ 2 h}  $ & 0 & 0 & 0 &0 \\
Lax & $\frac{1}{ \tau} $ &  0 & 0 & $  \frac{-1}{ 2 \tau} + \frac{c}{ 2 h } $ & $  \frac{-1}{2 \tau}  - \frac{c}{ 2 h} $ & 0 & 0 & 0 &0 \\
Lax-Wendroff & $\frac{1}{ \tau} $ & $  \frac{-1}{ \tau}  + \frac{c^2 \tau}{ h ^2} $ & 0 & $ \frac{( 1- \sigma ) c}{ 2 h} $&  $  \frac{-( 1+ \sigma ) c}{ 2 h } $ & 0 & 0 & 0 &0 \\
Crank-Nicolson & $ \frac{2}{ \tau}$ & $  \frac{2}{ \tau} $ & 0 & $ \frac{-c}{2h } $ & $ \frac{c}{ 2h } $ & $ \frac{c}{ 2h } $ & 0 & $ \frac{-c}{ 2h} $ & 0
\end{tabular}
\end{center}
\label{tab:scheme}
\end{table}%

\section{Caustics}
\label{sec:caust}

The solution of Eq. 
(\ref{transp}) is taken under the form:

\begin{equation}
\label{depl} u(x,t,k)=e^{i\,(k\,x-\omega \, t)}
\end{equation}

where $\omega = {\xi}_\omega+i \,{\eta}_\omega $ is the complex phase, and $k$ the real
wavenumber.
For dispersive waves, it is recalled that  the group velocity $V_g (k)$ is defined as 

\begin{equation}
\label{Vg} 
V_g (k) \equiv \frac{{\partial } \xi _\omega }{\partial {k}}  
\end{equation}

A caustic is defined as a focusing of different rays in a single location. The equivalent condition is that the group velocity exhibits an extremum, i.e. there exists at least one wave number $k_c$ such that

\begin{equation}
\label{eq:caustic} 
\frac{{\partial } V_g}{\partial {k}}  (k_c) =0
\end{equation}

The corresponding physical interpretation is that wave packets with characteristic wave numbers 
close to $k_c$ will pile-up after a finite time and will remain superimposed for a long time, 
resulting in the existence a region of high energy followed by a region with very low fluctuation level.

The linear continuous model Eq. (\ref{transp}) is not dispersive if the convection velocity $c$ is uniform, 
and therefore the exact solution does not exhibits caustics since the group velocity does not depends 
on $k$. The discrete solution associated with a given numerical scheme will admit spurious caustics, 
and therefore spurious local energy pile-up and local sudden growth of the error, if 
the discrete dispersion relation is such that the condition (\ref{eq:caustic}) is satisfied.
For a uniform scale-dependent convection velocity, such spurious caustics can exist in polychromatic 
solutions only, since they are associated to the superposition of wave packets with different 
characteristic wave numbers.

The general dispersion relation associated with the discrete scheme (\ref{eq:discrete-scheme}) is 

\begin{equation}
 \alpha \, e^{\,i\, \varphi } +\{\gamma \,e^{\,i \,\varphi }+\zeta\}\,e^{\,2\, i \,\omega \,\tau
   } +\{(\beta +e^{\,i \,\varphi }\,\delta )+\varepsilon \}e^{\,i \,\varphi
   }e^{\,i\,\omega \,\tau}\text
   +\eta\, e^{\,2 \,i \varphi }+\theta
=0
\end{equation}

 which is a non linear quadratic equation in $e^{i \omega \,\tau   }$ that can easily been solved. 

 The corresponding group velocity is given by:
 
\begin{equation}
\label{eq:Vg} 
V_g =h  \frac{{\partial} \xi _\omega }{\partial {\varphi}}
  \end{equation}
  
The numerical solution will therefore admits spurious caustics if

\begin{equation}
\frac{{\partial } V_g }{\partial {k}} = 
\frac{{\partial } V_g }{\partial {\varphi}} \, \frac{{\partial }
\varphi }{\partial {k}}=0
 \Longleftrightarrow
\frac{{\partial } V_g }{\partial {\varphi}}=0
          \end{equation}
          
The corresponding values of $\varphi$ and $k$ will be
respectively denoted
${\varphi}_c$ and $k_c$.

Spurious caustics are associated with characteristic lines given by 

\begin{equation}
 \frac{x}{t}=U_c
\end{equation}

 where
 
\begin{equation}
 U_c=V_g(\varphi_c)
\end{equation}

We now illustrate the caustic phenomenon considering  the two following sinusoidal wave packets:

\begin{equation}
\label{paquet} u_1=e^{\,-\alpha\,
(x-x_0^1-c\,t)^2}\,\text{Cos}\,[\,k_1\,(x-x_0^1-c\,t)\,] \,\,\, ,
\,\,\, u_2=e^{\,-\alpha\,
(x-x_0^2-c\,t)^2}\,\text{Cos}\,[\,k_2\,(x-x_0^2-c\,t)\,]
\end{equation}

where $\alpha>0$).
The two wave packets are initially centered  at $x_0^1$ and $x_0^2$, respectively.


If the solution obeys the linear advection law given by Eq. (\ref{transp}), the initial field is 
passively advected at speed $c$, while, 
if the advection speed is scale-dependent (as in numerical solutions), the two packets will travel at different speeds, leading to the rise
of discrepancies with the constant-speed solution. 
Another dispersive error is the shape-deformation phenomenon: due to numerical errors, the exact shape of the wave packets will not be exactly preserved.
This secondary effect will not be considered below, since it is not related to the existence of spurious caustics.

The spurious caustic will appear if the two wave packets happen get superimposed. During the cross-over, the $L_\infty$ norm of the error 
(defined as the difference between the constant-speed solution and the dispersive one) will exhibit a maximum. 
The  characteristic life time of the caustic depends directly on the difference 
between the advection speeds of the two wave packets and the wave packet widths.

Neglecting shape-deformation effects and assuming that the numerical scheme is non-dissipative, 
the numerical error $E$ is given by:

\begin{equation}
\begin{aligned}
E = \vert & \, e^{- \alpha  (x-x_0^1-c\,t )^2}
\text{Cos} [ k_1\,(x-x_0^1-ct)\,] -  e^{- \alpha
(x-x_0^1-t \,{V_1})^2} \text{Cos} [ k_1 (x-x_0^1-t {V_1})] \\
&+ e^{- \alpha (x-x_0^2-c\,t)^2} 
\text{Cos} [ k_2 (x-x_0^2-c t) ] - e^{ - \alpha (x-x_0^2-t \,{V_2})^2} \, \text{Cos}\,[\,k_2\,(x-x_0^2-t
\,{V_2})\,] \vert
\end{aligned}
\end{equation}

where $V_1$ and $V_2$ are the advection velocity of the two wave packets, respectively.

A simple analysis show that

\begin{equation}
\lim _{t \rightarrow + \infty} L_\infty ( E(t) ) = L_\infty ( u_1 (t=0)) , \quad \max _t  L_\infty ( E(t) ) = 2 L_\infty ( u_1 (t=0))
\end{equation}

\begin{figure}[htbp]
\centerline{\psfig{file=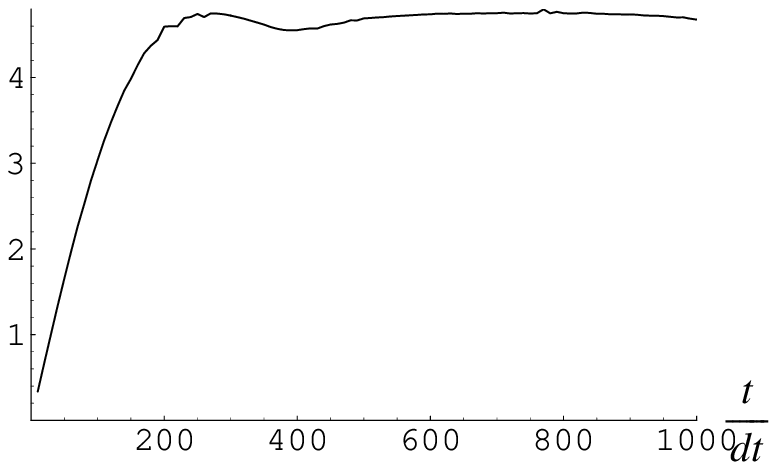,width=7cm} }
\centerline{\psfig{file=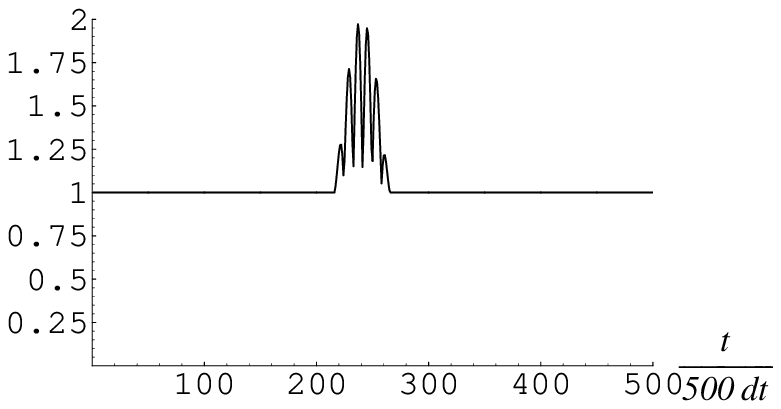,width=7cm}}
\caption{Time history of the numerical error for the two-wave packet problem (shape deformation and 
dissipative errors are neglected to emphasis the linear focusing phenomenon). 
Top: $L_1$ norm. Bottom: $L_\infty$ norm. 
Numerical parameters are $\alpha =1$, $h = 0.01$, $V_1 = 2.04$, $V_2 = 2.02$, 
corresponding to the properties of the Lax scheme for $\sigma =0.7$}
\label{fig-error}
\end{figure}

The time histories of the $L_1$ and $L_\infty$ norms of $E$ are displayed in Fig. \ref{fig-error}, 
showing the occurance of the caustic and the sudden growth of the $L_\infty$ error norm.

\section{Analysis of test numerical schemes}
\label{sec:analyze-test}

We analyzed  the properties of the numerical scheme displayed in section \ref{sec:scheme}.

For the Leapfrog scheme, the dispersion relation is:

\begin{equation}
1- e^{\,2\,i\,\omega  \,\tau}+ 2\,i\,\sigma \,\text{Sin}\,[\varphi ]
   \,e^{\,i\,\omega  \,\tau} =0
\end{equation}

from which it comes that

\begin{equation}
 e^{\eta _{\omega }\,\tau } =1 \Longrightarrow  \eta _{\omega } = 0
\end{equation}

\begin{equation}
{\,\xi _{\omega }\,\tau }=\,\text{ArcSin}\,[\sigma\,\varphi ]
\end{equation}

The group velocity can be expressed as
 
\begin{equation}
V_g = \frac{c \,
\text{Cos}\,[\varphi ]}{\sqrt{1-\sigma ^2 \,
\text{Sin}^2\,[\varphi ]}}
\end{equation}

leading to

\begin{equation}
\label{eq:leapfrog-grad}
\frac{{\partial } V_g }{\partial {\varphi}}= \frac{c \, h \text{Sin} [\varphi ]
}{  {\sqrt{1-\sigma ^2 \,  \text{Sin}^2\,[\varphi ]}} }  {\left(\frac{  \sigma ^2 \,  \text{Cos}^2\,[\varphi]} {\left(1-\sigma ^2 \, \text{Sin}^2\,[\varphi
   ]\right)}-1 \right)}
\end{equation}

A trivial root is   $\varphi=0 \, \, mod \, \,\pi$, which corresponds to

\begin{equation}
{U_c}= \pm \, c
\end{equation}

\begin{figure}[htbp]
\centerline{\psfig{file=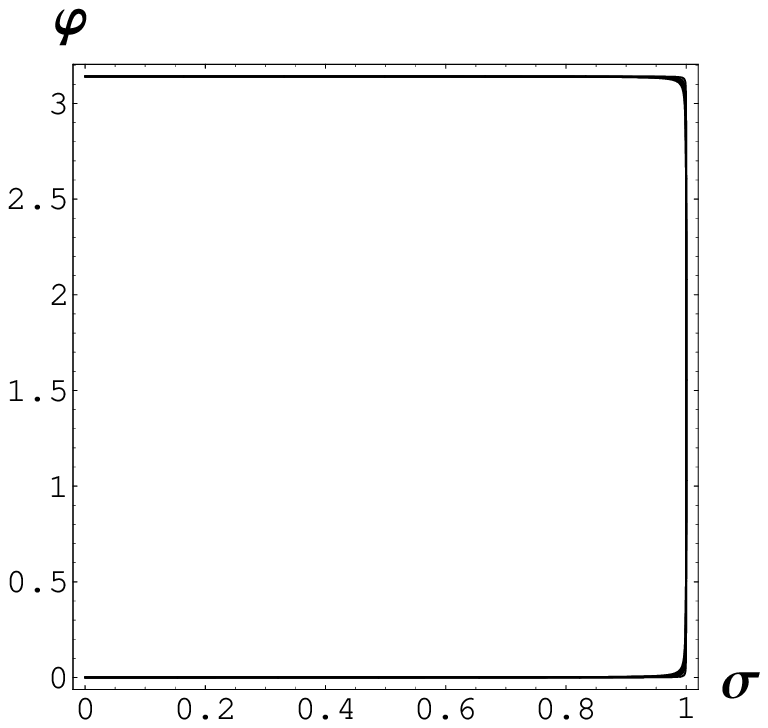,width=7cm}}
\caption{Roots of Eq. (\ref{eq:leapfrog-grad}), which correspond to possible occurance of spurious caustics for the Leapfrog scheme.}
\label{fig:roots-leapfrog}
\end{figure}

Zeros of  Eq. (\ref{eq:leapfrog-grad}) in the $(\sigma , \varphi )$ plane 
are displayed in Fig. \ref{fig:roots-leapfrog}.
  
Let us now consider  the Lax scheme.
The dispersion relation is:

\begin{equation}
\label{dispLax} 
2 \,e^{\,i \,\varphi }+\{(\sigma -1) \, e^{\,2\,\,
\varphi }-(\sigma +1)\}\,e^{\,i\,\omega  \,\tau} =0
\end{equation}

which yields:

\begin{equation}
\label{etaLax}
 e^{\,\eta _{\omega }\,\tau } = {\sqrt{{\text{Cos}}^2\,[\varphi
   ]+\,{\sigma} ^2\, {\text {Sin}}^2\,[\varphi ]}}
\end{equation}

\begin{equation}
{\xi _{\omega }\,\tau }  =  \text{ArcTan}\,\big[\,\sigma\,\text
{Tan}\,[\varphi ]\,\big]
\end{equation}

The group velocity is equal to

\begin{equation}
\label{VgLax} 
V_g =  \frac{c }{\text{Cos}^2[\varphi ]+\sigma ^2\,  \text{Sin}^2[\varphi ]}
\end{equation}

yielding

\begin{equation}
\label{eq:lax-grd}
\frac{{\partial } V_g }{\partial {\varphi}}= \frac{2 \,\tau^3
\left(\sigma ^2-1\right)\, \text{Sin}[\varphi
]\,\text{Cos}\,[\varphi
   ]}{ {\sigma }\, \left(\text{Cos}^2\,[\varphi ]+\sigma ^2 \,\text{Sin}^2\,[\varphi]\right)^2}
          \end{equation}
          
\begin{figure}[htbp]
\centerline{\psfig{file=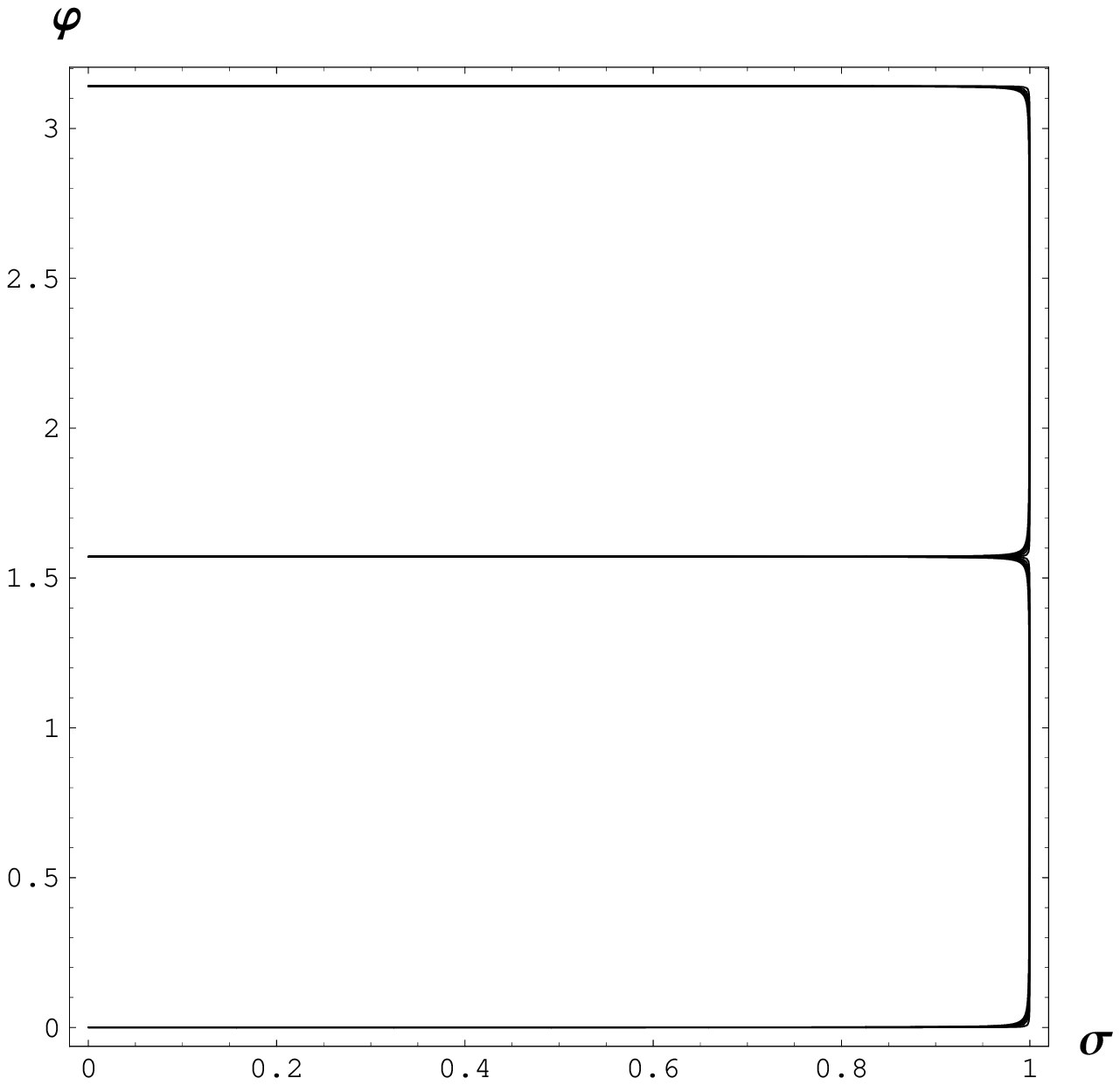,width=7cm}}
\caption{Roots of Eq.(\ref{eq:lax-grd}), which correspond to possible occurance of spurious caustics for the Lax scheme.}
\label{fig:roots-lax}
\end{figure}

Spurious caustics (see Fig. \ref{fig:roots-lax} ) can arise  for $\varphi=0 \, \, mod \, \,\pi$ and $\varphi=\frac {\pi}{2} \, \, mod \, \,\pi$.\\
The corresponding group velocities are respectively:

\begin{equation}
\label{UcLax} {U_c}^1  c, \quad   {U_c}^2  \frac{c }{\sigma ^2}
  \end{equation}

For the Lax-Wendroff scheme, the dispersion relation is:

\begin{equation}
\label{dispLaxWendroff} 
2 \,e^{i\, \,\varphi }+\sigma \,(\sigma -1)\,e^{\,2\,i\, \varphi } - \sigma \,
   (\sigma +1)+2  \,(\sigma ^2-1)\,e^{\,i \,\varphi }\,e^{\,i\,\omega \, \tau\,}=0
\end{equation}

which yields:

\begin{equation}
 e^{\,\eta _{\omega }\,\tau }  = 
 \sqrt{1-4\,(1-\sigma ^2)\,\sigma ^2 \,\text{Sin}^4[\frac{\varphi }{2}]}
\end{equation}

\begin{equation}
{\xi _{\omega }\,\tau }  =  \frac{\sigma \, \text{Sin}\,[\varphi ]}
{1-2\,\sigma ^2\,\text{ Sin}^2\,\left[\frac{\varphi }{2}\right]}
\end{equation}

The group velocity is given by 

\begin{equation}
V_g = c  \frac{ \left\{ \left(1-2 \sigma ^2 \, \text{Sin}^2\left[\frac{\varphi
   }{2}\right]\right)\text{Cos}\,[\varphi ]+\sigma ^2 \,\text{Sin}\,[\varphi ]\,
   \text{Sin}^2\left[\frac{\varphi }{2}\right]\right\}} {\left(1-2 \, \sigma
   ^2\,
   \text{Sin}^2\,\left[\frac{\varphi }{2}\right]\right)^2+\sigma ^2 \,\text{Sin}\,[\varphi
   ]^2}
\end{equation}

from which

\begin{equation}
\label{eq:laxw-grd}
\frac{{\partial } V_g }{\partial {\varphi}}=  - \, \frac{c^2 \, \tau
\, \left(-1+\sigma ^2\right) \left(-2+\sigma ^2+3 \, \sigma ^4-4
\,   \sigma ^4 \, \text{Cos}\,[\varphi ]+\sigma ^2 \,
\left(-1+\sigma ^2\, \right)\,
   \text{Cos}\,[2\,
   \varphi ]\right) \, \text{Sin}\,[\varphi ]}{2 \, \sigma \,  \left(\left(1-\sigma ^2+\sigma
   ^2 \, \text{Cos}\,[\varphi ]\right)^2+\sigma ^2 \, \text{Sin}^2\,[\varphi ]\right)^2}
          \end{equation}
          
\begin{figure}[htbp]
\centerline{\psfig{file=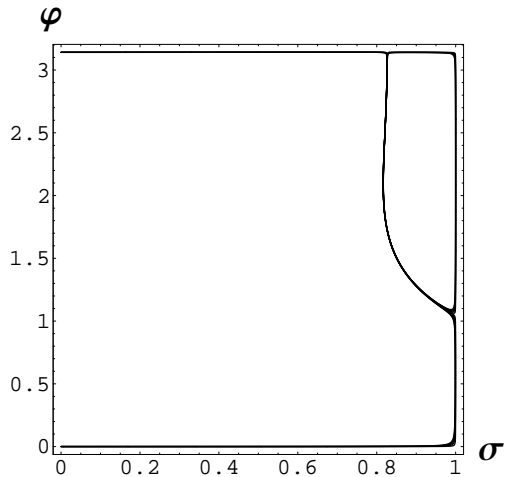,width=7cm}}
\caption{Roots of (\ref{eq:laxw-grd}), which correspond to possible occurance of spurious caustics for the Lax-Wendroff scheme.}
\label{fig:roots-laxw}
\end{figure}

Roots of  Eq. (\ref{eq:laxw-grd}) are shown in Fig. \ref{fig:roots-laxw}.

We finally consider the Crank-Nicolson scheme, whose dispersion relation is:

\begin{equation}
\label{dispCrank} 
\sigma ^2 \, (\text{ Cos }\,[\varphi ]-1)
\,\text{ Cos }[\frac{\omega  \, \tau }{2}]+i \,c  \,\tau
   \text{ Sin }[\frac{\omega  \,\tau}{2}]=0
\end{equation}

which yields:

\begin{equation}
e^{\eta _{\omega }\,\tau }  =  \frac{\pm (c \,\tau-2\,\sigma
^2 \, \text{Sin}^2\,[\frac{\varphi
   }{2}])}{c \,\tau+2\,\sigma ^2 \,
   \text{Sin}^2 \,[\frac{\varphi }{2}]}
\end{equation}

\begin{equation}
\xi _{\omega }\,\tau  =  0 \,\, mod \,\, \pi
\end{equation}

The group velocity being constant, no spurious caustics can arise with this scheme.

\section{Concluding remarks}

The existence of spurious numerical caustics in linear advection schemes has been proved.
This linear dispersive phenomenon gives rise to a sudden growth of the $L_\infty$ norm of the error, 
which corresponds to a local focusing of the numerical error in both space and time.
In the present analysis, spurious caustics have been shown to arise in polychromatic solutions.

The energy of the caustic phenomenon depends on the number of spectral modes that will get superimposed at the same time.
As a consequence, the spurious error pile-up will be more pronounced  in simulations with very small wave-number increments.

It has been shown that most popular existing schemes allow the existence of spurious caustics, while some schemes are caustic-free, like the Crank-Nicolson scheme.

\end{document}